# ABOUT A CLASS OF TRUNCATED OPERATORS REPRESENTED VIA SCHAUDER BASES WITH APPLICATIONS

M. De la Sen

Instituto de Investigacion y Desarrollo de Procesos. Universidad del Pais Vasco

Campus of Leioa (Bizkaia) - P.O. Box 644- Bilbao, 48080- Bilbao. SPAIN

**Abstract**. This research investigates the relations between the fixed points of asymptotically contractive operators defined on complete liner metric spaces, in general, and in separable Hilbert spaces in particular and their obtained truncations of the expansions of such operators by using Schauder bases.

## 1. Introduction

Separable Hilbert spaces are very important for the study of certain relevant applied problems in Control Theory and Signal Theory, [1]. It is well-known that operators on separable Hilbert spaces can be represented in closed form through expansions using either orthogonal or orthonormal bases. Another representative property associated with such operators is that they admit truncations using a finite number of the members of the chosen basis [1-2]. On the other hand, attention has been devoted to the properties concerning the existence and uniqueness of fixed points of mappings in metric, Banach and more general spaces. The relevance of the theory in properties has been studied concerning theoretical aspects and their applications such as the existence and uniqueness of solutions in differential, difference and hybrid equations as well as in continuous-time, discrete -time, digital and hybrid dynamic systems stability theory in the above problems [5, 9-11, 17], [21-26]. Important issues related to stability studies are the existence/ uniqueness of fixed points and best proximity points, as well as the boundedness of the iterated sequences which are built via " ad hoc" maps together with the convergence of such iterated sequences to limit points. See, for instance, [3-17] and references there in. Importance has also been paid to the investigation of the existence and uniqueness of coupled and common fixed points and best proximity points for several mappings and related properties as well as to the relevance of fixed points in the context of variational inequalities [5-7, 8, 13, 20]. This paper is devoted to the study of some fixed point properties of self-mappings and operators defined on complete linear metric spaces, in general, and in separable Hilbert spaces and the comparison to those properties being associated with their finite-dimensional truncated counterparts by taking a finite number of the members of the Schauder basis. The major highlights of the paper are the ability of approximately maintaining the asymptotic properties of the distances as well as either the convergence or the approximate convergence of sequences built through the truncated maps to defined neighborhoods of the fixed points to which the sequences built via the untruncated maps converge.

## 2. Preliminaries

Consider a separable complete linear metric space $(X,d)$ and a sequence $\{e_n\} \subset X$ being a Schauder basis (or simply a basis of the space $X$), that is, there is a unique representation $x = \sum_{i=1}^{\infty} x_i e_i$; $\forall x \in X$ where the coefficients $x_i = x_i(x)$ are scalars. Let us define $[x]_m := P_m x = \sum_{i=1}^{m} x_i e_i$ in the subspace $X_m$ of $X$ where $P_m : X \to X_m$ is a linear operator for any $m \in N$ and note that $x = [x]_m + \sum_{i=m+1}^{\infty} x_i e_i$.

It is known that the sequence $\{P_m(x)\}$ is uniformly bounded by some constant $C$. If $C=1$ the Schauder basis is monotone and then the elements of $\{P_m(x)\}$ are the Schauder basis projections. Schauder basis can also be defined on Banach spaces $(X,\|\ \|)$ such that the linear space $X$ (defined on some field $F$) is separable, the $x_i = x_i(x) \in F$ are unique for each $x \in X$ and the convergence is understood with respect to the norm topology, i.e. $\left\| x - \sum_{i=1}^n x_i e_i \right\|_X \to 0$ as $n \to \infty$. It is also known that all the operators $P_m: X \to X_m$ are equicontinuous since for each neighborhood of zero $V_m \subset X_m$, there is a neighborhood of zero $U \subset X$ such that $P_m U \subset V$, [2]. Note that in order that a linearly independent set $\{e_n\}$ be a Schauder basis of a linear metric space $(X,d)$, such a set has to be dense in $X$ and the operators $P_m x = [x]_m$ are equicontinuous in the space $lin\{e_n\}$. If $T$ is a mapping of $X$ into itself then one has, in general, for arbitrary $x \in X$

$$T[x]_m = T\left(\sum_{i=1}^m x_i e_i\right) \neq [Tx]_m = \sum_{i=1}^m (Tx)_i e_i$$

Since, $T(x_i e_i) \neq (Tx)_i e_i$, in general. However, if the equality holds for some $x \in X$ then it holds for any functional powers (or iterated functions) of $T$, i.e. any repeated composition of any iteration order of $T$ with itself for $x \in X$ as it is proven in the subsequent result.

**Proposition 2.1**. Let $(X,d)$ be a complete linear metric space and let $\{e_n\} \subset X$ be a Schauder basis. Assume that $T(x_i e_i) = (Tx)_i e_i$ for some $X \ni x = \sum_{i=1}^\infty x_i e_i$ and some $m \in N$; $\forall i \in \overline{m} := \{1, 2, \ldots, m\}$. Then, $(T^n x)_i e_i = T^n(x_i e_i)$; $\forall i \in \overline{m}$, $\forall n \in N$. If, in addition, $T^n : X \to X$ is linear then $T^n [x]_m = [T^n x]_m$. □

**Theorem 2.2**. Let $\hat{X}$ be the linear space of all sequences $\{t_i\}$ of scalars such that the series $\sum_{i=1}^\infty t_i e_i$ is convergent with a norm on $\hat{X}$ being defined by $\|y\|_{\hat{X}} := \sup_{m \in N} \left| \sum_{i=1}^m t_i e_i \right|$, where $\{e_n\} \subset X$ is a Scahauder basis and let $T: X \to X$ be an asymptotically nonexpansive self-mapping and $d: X \times X \to R_{0+}$ be the $\|.\|_{\hat{X}}$-induced metric. Then, the subsequent properties hold:

**(i)** $\lim_{m \to \infty} \left| \left( (T^{n+k} x)_m - (T^n x)_m \right) e_m \right| = 0$; $\forall n \in N$, $\forall x \in X$

**(ii)** For any given $\varepsilon \in R_+$, there are $n_0 = n_0(\varepsilon,x) \in N$ and $m_0 = m_0(\varepsilon,x) \in N$ such that

$$d\left( [T^{n+k} x]_m, [T^n x]_m \right) \leq (1+\varepsilon) d(x, T^k x) + \varepsilon \ ; \ \forall m, n(\geq n_0) \in N, \forall x \in X$$

and

$$\limsup_{n,m \to \infty} \left( d\left( [T^{n+k} x]_m, [T^n x]_m \right) - d(x, T^k x) \right) \leq 0; \ \forall k \in N, \forall x \in X$$

**(iii)** If $T: X \to X$ is asymptotically contractive then

$$d\left( [T^{n+k} x]_m, [T^{n+k} x]_m \right) \to 0 \text{ as } n, m \to \infty; \ \forall k \in N, \forall x \in X$$

**(iv)** Assume that $T: X \to X$ is linear and $T(x_n e_n) = (Tx)_n e_n$. Then, Properties (i)-(iii) also hold "mutatis-mutandis" by replacing $[T^{n+k} x]_m \to T^{n+k} [x]_m$. □

**Proposition 2.3.** Let $(X, d)$ be a complete linear metric space, let $\{e_n\} \subset X$ be a Schauder basis and assume that $T: X \to X$ is asymptotically contractive. If the metric subspace $(X_m, d_{X_m})$ of $(X, d)$ has a fixed point of $T: X | X_m \to X_m$ then it is necessarily unique and coincident with that of $X$

*Proof:* Given $x = \sum_{i=1}^{\infty} x_i e_i$ in $X$, take $[x]_m := P_m x = \sum_{i=1}^{m} x_i e_i$ in $X_m \subseteq X$. Since $T: X \to X$ is asymptotically contractive there is a real constant $K \in [0,1)$ and a sequence of finite natural numbers $\{p_n\}$ such that $d(T^n x, T^n y) \leq K^{n-p_n} d(T^{p_n} x, T^{p_n} y); \forall x, y \in X$ with $p_n = p_n(n)$ and $\overline{p} := \limsup_{n \to \infty} p_n < +\infty$. Since the sequence $\{p_n\}$ is bounded, there is a finite nonnegative real constant $M = M(x, y) \geq \sup_{n \in N} d(T^{p_n} x, T^{p_n} y)$ leading to the existence of the limit

$$\lim_{n \to \infty} d(T^n x, T^n y) = \lim_{n \to \infty} d(T^n x, T^{n+1} x) = d\left(\lim_{n \to \infty} (T^n x), T^{n+1} x\right)$$

$$= d\left(\lim_{n \to \infty} (T^n x), T \lim_{n \to \infty} (T^n x)\right) = d\left(\lim_{n \to \infty} (T^n x), \lim_{n \to \infty} (T^n y)\right) = 0$$

; $\forall x, y \in X$. Then, $\{T^n x\} \to z = Tz (\in X)$ as $n \to \infty$ which is the unique fixed point of $T: X \to X$ in $X$; to which all the sequences $\{T^n x\}$ converge; $\forall x \in X$. Assume that there is some $m \in N$, such that $\{T^n y\} \to [T^n z_1]_m = [z_1]_m (\in X_m)$ for some $y \in X$, and that $z_1 \neq z$. Since $\{e_n\} \subset X$ is a Schauder basis and $X_m \subseteq X$, $z_1$ has the representation $z_1 = \sum_{i=1}^{m} z_{1i} e_i = \sum_{i=1}^{\infty} z_{1i} e_i$ in both $X_m$ and $X$ so that $z_{1i} = 0$ for $i > m$. Thus, $[T^n z_1]_m = T^n z_1 = z_1 = \sum_{i=1}^{m} z_{1i} e_i = \sum_{i=1}^{\infty} z_{1i} e_i$ so that, by using the contractive condition, one gets:

$$0 \leftarrow d(z_1, T^n y) = d(z, T^n y) = d(T^n z, T^n y) = d(T^n z_1, T^n y)$$

$$= d\left([T^n z_1]_m, T^n y\right) \leq K^{n-p_n} d(T^{p_n} z_1, T^{p_n} y) \to 0; \forall y \in X \text{ as } n \to \infty$$

so that $z = z_1 \in X_m$ is the unique fixed point of $T: X \to X$ in $X$. □

The key point in the above result, rather than the asymptotic contractiveness of the self-mapping, which guarantees the result is the fact that $T: X \to X$ has a unique fixed point $z$. Then, if $T: X | X_m \to X | X_m$ has a fixed point then such a point is $z$. Note that last result does not hold for a nonexpansive mapping with more than one fixed point.

**Example 2.4.** Consider the functional differential equation subject to a piecewise continuous time-varying point delay $\lambda(t) \in (0, h]$; $\forall t \in \mathbf{R}_{0+}$:

$$\dot{y}(t) = a\, y(t) + a_0(t)\, y(t - \lambda(t))$$

; $\forall t \in \mathbf{R}_{0+}$ where $a_0 : \mathbf{R}_{0+} \to \mathbf{R}$ is bounded piecewise-continuous under initial conditions $\varphi(t) \equiv x(t)$ for $t \in [-h, 0]$ with $\varphi : \mathbf{R}_{0+} \to \mathbf{R}$ being an absolutely continuous function, $a < 0$ and $\sup_{t \in \mathbf{R}_{0+}} |a_0(t)| < |a|$. Assume also that there is some positive real constant $T$ such that $\limsup_{n \to \infty} (T + \lambda[(n+1)T] - \lambda(nT)) \leq 0$ implying that $\{\lambda(nT)\} \to 0$. Consider the Hilbert space $(X, |\;|)$, which is also a complete linear metric space, of square-integrable solution sequences on $\mathbf{R}_{0+}$ with respect to the distance in-between points defined by the absolute value of the difference. Let $\{y_n \equiv y(nT)\}$ a sequence of samples of the unique solution for the inter-sample period $T > h$ for any $n \in \mathbf{N}$. It turns out that for any given real $\sigma \in [0, T)$ and any $n(\geq n_0) \in \mathbf{N}$:

$$\left| y[(n+1)T + \sigma - \lambda((n+1)T)] \right| = \left| e^{-|a|(T+\sigma-\lambda((n+1)T))} \left( y(nT + \sigma - \lambda((n+1)T)) + \int_0^T e^{|a|\tau} a_0(\tau) y(nT + \sigma + \tau - \lambda(\tau)) d\tau \right) \right|$$

$$\sup_{\sigma \in [0,T)} \left| y[(n+1)T + \sigma - \lambda((n+1)T)] \right| = \sup_{\sigma \in [0,T)} \left| e^{-|a|(T+\sigma-h)} \left( y(nT + \sigma - h) + \int_0^T e^{|a|\tau} a_0(\tau) y(nT + \sigma + \tau - \lambda(\tau)) d\tau \right) \right|$$

$$\leq \left( \sup_{\sigma \in [0,T)} e^{-|a|(\sigma - h)} \right) \sup_{\sigma \in [0,T)} |y(nT + \sigma - \lambda(nT))|$$

Thus, $y(t) \to 0$ as $t \to \infty$ so that $y^* = 0$ is the unique equilibrium point of the differential equation and the unique fixed point of the self-mapping $T : \mathbf{N} \times X \to \mathbf{N} \times X$ mapping the each solution strip $y(nT + \sigma - \lambda(nT))$ into the next one $y((n+1)T + \sigma - \lambda((n+1)T))$ for each $k \in \mathbf{N}$ and each given real constant $\sigma > h$ and the unique fixed point of the self-mapping that generates the unique solution trajectory for any given function of initial conditions. Let $\{e_n\} \subset \mathbf{N} \times X$ be a Schauder basis, where $e_1 = e^{anT}$ and note that $T : \mathbf{N} \times X \to \mathbf{N} \times X$ is asymptotically contractive. Let $Y_\sigma \equiv \{y(nT + \sigma)\}$. Then, $y = [y]_1 + \sum_{i=2}^{\infty} y_i e_i$; $\forall y \in Y_\sigma$ where $[y]_1 \in X_1$ is the truncation of $y$ for $m = 1$ if $y_1 = \varphi(0)$. Note that Proposition 2.3 follows with $[y]_1 (nT + \sigma) = e^{-|a|(nT+\sigma)} \varphi(0) \to y^* = 0$ as $n \to \infty$, which is the solution of the delay-free differential equation $\dot{x}(t) = a\, x(t)$ for $x(0) = y(0) = \varphi(0)$, $y(nT) \to y^* = 0$ as $n \to \infty$, and

$$\lim_{n \to \infty} \left( \left( y - [y]_1 - \sum_{i=2}^{\infty} y_i e_i \right) (nT + \sigma) \right) = \lim_{n \to \infty} \left( (y(nT + \sigma) - [y(nT + \sigma)]_1) - \sum_{i=2}^{\infty} y_i (nT + \sigma) e_i (nT + \sigma) \right)$$

$$= \lim_{n \to \infty} \left( \int_{nT}^{nT + \sigma} e^{-|a|(nT + \sigma - \tau)} a_0(\tau) y(\tau - \lambda(\tau)) d\tau \right) = 0 \qquad \square$$

## 3. Some extensions

The next result offers some parallel properties to Theorem 2.2 for the case of linear metric spaces endowed with a homogeneous and translation-invariant metric.

**Theorem 3.1**. Let $(X, d)$ be a linear metric space endowed with a homogeneous and translation-invariant metric $d : X \times X \to \mathbf{R}_{0+}$ and let $\{e_n\} \subset X$ be a Schauder basis of linearly independent elements of $X$. The following properties hold if $T : X \to X$ is an asymptotically contractive self-mapping:

(i) $\lim\limits_{n,m \to \infty} d\left([T^n x]_m, [T^n y]_m\right) = 0$; $\forall x, y \in X$.

(ii) Assume that $(X, d)$ is a complete linear metric space and $D$ is any given bounded domain of $X$ then, for any given $\varepsilon \in \mathbf{R}_+$, there is a natural number $\overline{m} = \overline{m}(\varepsilon, D)$ such that $d\left([T^n x]_m, [T^n y]_m\right) \leq \varepsilon$, $d(T^n x, T^n y) \leq \varepsilon$ and $d\left([T^n x]_m, T^n x\right) \leq 2\varepsilon$; $\forall n, m \geq \overline{m}$, $\forall x, y \in D$.

(iii) If $(X, d)$ is complete and $z \in X$ is the (unique) fixed point of $T : X \to X$ then for the given $\varepsilon \in \mathbf{R}_+$ in Property (ii), there is a neigborhood $B(z, 4\varepsilon) = \{r \in X : d(r, z) < 4\varepsilon\}$ of $z$ of radius $4\varepsilon$ such that $\left\{[T^{n+1} x]_m\right\} \in cl\, B(z, 4\varepsilon)$; $\forall n, m \geq \overline{n}$, $\forall x \in D$ and some $\overline{n}(\geq \overline{m}) \in \mathbf{N}$. As a result, $\left\{[T^n x]_m\right\} \to z$ as $n, m \to \infty$; $\forall x \in D$.

(iv) Assume that $T : X \to X$ is linear and $T(x_n e_n) = (Tx)_n e_n$. Then, Properties (i)-(iii) also hold "mutatis-mutandis" by replacing $[T^{n+k} x]_m \to T^{n+k} [x]_m$.

*Proof*: Since the linear metric space $(X, d)$ is endowed with an homogeneous and translation-invariant metric $d : X \times X \to \mathbf{R}_{0+}$ then the linear metric space is also a normed space $(X, \|\,\|) \equiv (X, d)$ and, since $T : X \to X$ is asymptotically contractive, $d(T^n x, T^n y) \leq K^{n-p_n} d(T^{p_n} x, T^{p_n} y)$; $\forall x, y \in X$ with $p_n = p_n(n)$ and $\overline{p} := \limsup\limits_{n \to \infty} p_n < +\infty$. Since the sequence $\{p_n\}$ is bounded, there is a finite nonnegative real constant $M = M(x, y) \geq \sup\limits_{n \in \mathbf{N}} d(T^{p_n} x, T^{p_n} y)$, one has for some real constant $K \in [0, 1)$

$$d\left([T^n x]_m, [T^n y]_m\right) \leq d\left([T^n x]_m, T^n x\right) + d(T^n x, T^n y) + d\left([T^n y]_m, T^n y\right)$$

$$\leq K^{n-p_n} M(x, y) + \|T^n\| \left(\sum_{i=m+1}^{\infty} (|x_i| + |y_i|) \|e_i\|\right); \forall x, y \in X$$

Since the Schauder basis $\{e_n\}$ of $X$ consists of linearly independent elements of $X$ then any linear combination of the elements $e_n$ is dense in the space $X$. Furthermore, since the metric is homogeneous and translation-invariant it has an induced norm which is a $p$-homogeneous norm for $p = 1$ (abbreviated simply as an homogeneous norm). Thus, from the above to joint properties, there is a nonnegative real number $C = C(x,m)$, dependent on $x$ and $m$, such that $\left\|\sum_{i=1}^{m} x_i e_i\right\| \leq C(x,m) \left\|\sum_{i=1}^{\infty} x_i e_i\right\|$; $\forall x \in X$. Thus,

$$\|T^n\|\left\|\sum_{i=m+1}^{\infty} x_i e_i\right\| \leq \|T^n\| C_1(x,m)\left(\left\|\sum_{i=1}^{m} x_i e_i\right\| + \left\|\sum_{i=1}^{\infty} x_i e_i\right\|\right) \leq \|T^n\| C_1(x,m)(1+C(x,m))\left\|\sum_{i=1}^{\infty} x_i e_i\right\|; \forall x \in X$$

since by construction, there is a nonnegative real number $C_1 = C_1(x,m)$, dependent on $x$ and $m$, such that such that the function $C_1(x,m)$ satisfies $\lim_{m \to 0} C_1(x,m) = 1$ and $\lim_{m \to \infty} C_1(x,m) = 0$, $\forall x \in X$. One gets that the limit below exists:

$$\lim_{n \to \infty} d(T^n x, T^n y) = \lim_{n,m \to \infty} d([T^n x]_m, [T^n y]_m) = 0 ; \forall x,y \in X$$

and Property (i) is proven. On the other hand, one gets:

$$0 \leftarrow d([T^n x]_m, [T^n y]_m) \leq K^{n-p_n} M(x,y) + \|T^n\| C_1(x,m)(1+C(x,m))\|x\|$$
$$+ \|T^n\| C_1(y,m)(1+C(y,m))\|y\| \to 0 ; \forall x,y \in X$$

as $n,m \to \infty$, $\forall x,y \in X$. Since $K^{n-p_n} M(x,y) \to 0$ as $n \to \infty$ and $\{g(x,m)\} \to 0$; $\forall x,y \in X$ where $g(x,m) := \sum_{i=m+1}^{\infty} (|x_i| + |y_i|)\|e_i\|$ then $\{g(x,m)\}$ is bounded, and for any nonzero $x \in X$ and any given $\varepsilon \in R_+$, there are $m_1 = m_1(\varepsilon, x, y) \in N$, $m_2 = m_2(\varepsilon, x) \in N$ and $m_3 = m_3(\varepsilon, y) \in N$ such that for any $n,m (\geq \max(m_0, m_1, m_2)) \in N$, one gets:

$$K^{n-p_n} M(x,y) \leq \varepsilon/3, \quad \|T^n\| C_1(x,m)(1+C(x,m))\|x\| \leq \varepsilon/3, \quad \|T^n\| C_1(y,m)(1+C(y,m))\|x\| \leq \varepsilon/3$$

For any given $\varepsilon \in R_+$, there is $N \ni \overline{m}' = m(\varepsilon, D) := \sup_{x,y \in D} \max(m_0(\varepsilon, x, y), m_1(\varepsilon, x, y), m_2(\varepsilon, x, y))$, where $D$ is any bounded domain of $X$, one has that $d([T^n x]_m, [T^n y]_m) \leq \varepsilon$; $\forall n,m (\geq \overline{m}') \in N$, $\forall x,y \in X$ and Property (ii) follows. By taking $m = \infty$, it follows that there is a $n_0' = n_0'(\varepsilon, x, y) \in N$ such that $d(T^n x, T^n y) \leq \varepsilon$; $\forall n (\geq n_0') \in N$. Thus, there is $\overline{m} \geq \max(\overline{m}', n_0')$, such that $\max\left[d([T^n x]_m, [T^n y]_m), d(T^n x, T^n y)\right] \leq \varepsilon$; $\forall n,m (\geq \overline{m}) \in N$, $\forall x,y \in D$. On the other hand, note that:

$$\left|d([T^{n+1} x]_m, T^{n+1} x) - d([T^n x]_m, T^n x)\right| \leq d([T^{n+1} x]_m, [T^n x]_m) + d(T^n x, T^{n+1} x) \leq 2\varepsilon$$

; $\forall n,m (\geq \overline{m}) \in N$, $\forall x \in D$. Property (ii) has been proven. Since $T : X \to X$ is asymptotically contractive, it has a unique fixed point $z$ in $X$ such that $\{T^n x\} \to z$ with $\{T^n x\}$ being a Cauchy sequence since any convergent sequence in a complete metric space is a Cauchy sequence. As a result, for, there is a $\bar{n} = \bar{n}(\varepsilon, \bar{m}, D)(\geq \bar{m}) \in N$ such that $max\left( d(T^n x, z), d\left( [T^n x]_m, [T^{n+1} x]_m \right) \right) \leq \varepsilon$ and

$d\left( T^n x, [T^n x]_m \right) \leq 2\varepsilon$; $\forall n, m > \bar{n}$, $\forall x \in D$ and then the proof of Property (iii) follows since

$$d\left( [T^{n+1} x]_m, z \right) \leq d\left( [T^n x]_m, [T^{n+1} x]_m \right) + d\left( [T^n x]_m, T^n x \right) + d\left( T^n x, z \right) \leq 4\varepsilon$$

Property (iv) is a direct consequence of Properties [(i)-(iii)] by using Proposition 2.1. □

**Theorem 3.2.** Let $(X, d)$ be a linear metric space endowed with a homogeneous and translation-invariant metric $d: X \times X \to \mathbf{R}_{0+}$ and let $\{e_n\} \subset X$ be a Schauder basis of linearly independent elements of $X$. Assume that $T: X \to X$ is a contractive self-mapping with contractive constant $K \in [0, 1)$. Then, the following properties hold:

$$d\left( [T^n x]_m, [T^n y]_m \right) \leq K^n \left( 1 + \sup_{i \leq n (\in N)} \left[ C_{1m}(T^i x, T^i y)(1 + C_m(T^i x, T^i y)) \right] \right) d(x, y); \forall x, y \in X$$

$$\limsup_{n,m \to \infty} \left[ d\left( [T^n x]_m, [T^n y]_m \right) - K d\left( [T^{n-1} x]_m, [T^{n-1} y]_m \right) \right] \leq 0; \forall x, y \in X$$

$$\limsup_{n,m \to \infty} \left[ d\left( [T^n x]_m, [T^{n+1} x]_m \right) - K d\left( [T^{n-1} x]_m, [T^n x]_m \right) \right] \leq 0; \forall x \in X$$

$$\exists \lim_{n,m \to \infty} d\left( [T^n x]_m, [T^{n+1} x]_m \right) = 0; \forall x \in X$$

$$\limsup_{n \to \infty} \left[ d\left( [T^n x]_m, [T^n y]_m \right) - d\left( [x]_m, [y]_m \right) - \frac{\bar{K}_{mE} \bar{\varepsilon}_{mE}}{1 - K} \right] \leq 0; \forall m \in N, \forall x, y \in E$$

where $E$ is any bounded subset of $X$ and $\bar{K}_{mE} \bar{\varepsilon}_{mE}$ are nonnegative real constants which depend on $m$ and $E$. If $(X, d)$ is complete and then $[T^n x]_m \to z$ as $n, m \to \infty$; $\forall x \in X$ where $z$ is the unique fixed point of $T: X \to X$. □